\pdfoutput=1
\documentclass[10pt,a4paper]{article}
\usepackage[T1]{fontenc}
\usepackage[left=2cm, right=2cm]{geometry}
\usepackage{amsthm}
\usepackage{amsmath}
\usepackage{amsxtra}
\usepackage{amstext}
\usepackage{amssymb}
\usepackage{hyperref}
\theoremstyle{definition}
\newtheorem{defn}{Definition}[section]
\newtheorem{rem}{Remark}[section]
\theoremstyle{plain}
\newtheorem{theo}{Theorem}[section]
\newtheorem{lem}{Lemma}[section]

\numberwithin{equation}{section}
\begin{document}
	%Shortened Title : Hemi-Slant of lcK as Warped Products
	\title{\textbf{Hemi-Slant Submanifolds of lcK Manifolds as Warped Products}\footnote{\textbf{Keywords and phrases:} warped product submanifolds, locally conformal Kähler manifolds (lcK), hemi-slant submanifolds\par\textbf{2020 AMS Subject Classification:} 53C15, 53C40, 53C42}}
	\author{\textbf{Umar Mohd Khan \& Viqar Azam Khan}\\Department of Mathematics\\Aligarh Muslim University\\Aligarh-202002, India\\Email: umar.007.morpheus@gmail.com, viqarster@gmail.com}
	\date{}
	\maketitle
	\begin{abstract}
		We study immersions of a hemi-slant submanifold of lcK manifolds as a warped product with the leaves of the holomorphic (respectively slant) distribution warped and establish characterisation theorems and estimations for the squared length of the second fundamental form in both cases.
	\end{abstract}
	\section{Introduction}
	Vaisman introduced locally conformal Kähler(lcK) manifolds as a generalisation of Kähler manifolds \cite{Vaisman-lcK1,Vaisman-lcK2,Vaisman-lcK3,Vaisman-lcK4,Vaisman-lcK5,Vaisman-lcK6,Vaisman-lcK7}. An lcK manifold is a Hermitian manifold that can be written as the union of Kähler manifolds such that the lcK metric is locally conformal to these Kähler metrics. LcK manifolds are characterised by the existence of a globally defined closed 1-form $\omega$, called the \textit{Lee form}, such that the fundamental 2-form of the lcK metric satisfies $d\Omega=\Omega\wedge\omega$. The Lee form and its associated Lee Vector field play an important part in the geometry of lcK manifolds.\par
	From an extrinsic geometric standpoint, holomorphic and totally real submanifolds are important objects of study in the setting of almost Hermitian manifolds. Bejancu \cite{Bejancu-CR-I,Bejancu-CR-II} defined CR submanifolds as a generalisation of holomorphic and totally real submanifolds which were further studied by Chen \cite{Chen-CR-I,Chen-CR-II}. Later, Chen \cite{Chen-SlantImm,Chen-SlantSubmanifolds} extended the class of holomorphic and totally real submanifolds by introducing the notion of slant submanifolds. The concept was further generalised to pointwise slant submanifolds \cite{Chen-PointwiseSlant} by the same author. The study of CR submanifolds and slant submanifolds was later generalised by several authors to semi-slant submanifolds, hemi-slant submanifolds(also called pseudo-slant submanifolds) and bi-slant submanifolds, in various ambient manifolds.\par
	Papaghiuc \cite{Papaghiuc-SemiSlant} studied semi-slant submanifolds in almost Hermitian manifolds. Cabrerizo et al. \cite{Cabrerizo-SemiSlant-Sasakian,Cabrerizo-Slant-Sasakian} studied semi-slant submanifolds in Sasakian manifolds. Slant and semi-slant submanifolds in almost product Riemannian manifolds were studied in \cite{Atceken-Slant-RiemProd,Li_Liu-SemiSlant-LocProd,Sahin-Slant-AlmProdRiem}. Hemi-slant submanifolds were also studied in nearly Kenmotsu manifolds \cite{Atceken-HemiSlant-Kenmotsu}, LCS-manifolds \cite{Atceken-Slant_HemiSlant-LCS} and locally product Riemannian manifolds \cite{Tastan_Ozemdir-HemiSlant-LocProdRiem}.\par
	Bishop and O'Neill \cite{WarpedProd-Def} while studying examples of manifolds with negative sectional curvature, defined warped product manifolds by homothetically warping the product metric on a product manifold. Warped products are a natural generalisation of Riemannian products and they have found extensive applications in relativity. Most notably the Schwarzschild metric describing the gravitational field outside a spherical mass under certain assumptions and the Robertsen Walker metric (FLRW metric) are examples of warped product metrics. A natural example of warped product manifolds are surfaces of revolution. Hiepko \cite{Hiepko-WarpedProdCharacterisation} gave a characterisation for a Riemannian manifold to be the warped product of its submanifolds, generalising the deRham decomposition theorem for product manifolds. Later on Nölker \cite{Nolker-WarpedExtrinsic} and Chen \cite{Chen-TwistedProd,Chen-WarpedProd1,Chen-WarpedProd2} initiated the study of extrinsic geometry of warped product manifolds.\par
	Chen \cite{Chen-WarpedProdCR-CRWarpedProd-I,Chen-WarpedProdCR-CRWarpedProd-II} initiated the study of CR submanifolds immersed as warped products in Kähler manifolds. Bonanzinga and Matsumoto \cite{Bonanzinga_Matsumoto-WarpedProdCR-lcK,Matsumoto_Bonanzinga-DoublyWarpedProdCR-lcKSpaceForm,Matsumoto_Bonanzinga-DoublyWarpedProdCR-lcKSpaceForm-II} continued the study in the setting of lcK manifolds. Nargis Jamal et al. \cite{NargisJamal_KAK_VAK-GenericWarped-lcK} studied Generic warped products in lcK manifolds. Further studies of semi-slant and hemi-slant submanifolds of lcK manifolds were carried out in \cite{Alghamdi_Uddin-SemiSlantWarped-lcK,Tastan_Sibel-HemiSlant-lcK,Tastan_Tripathi-SemiSlant-lcK}. Generic submanifolds, CR-submanifolds and semi-slant submanifolds immeresd as warped products in lcK manifolds were studied by \cite{NargisJamal_KAK_VAK-GenericWarped-lcK,Alghamdi_Uddin-SemiSlantWarped-lcK}.\par
	We continue the study by considering hemi-slant submanifolds in an lcK manifold. In particular we give characterisation theorems and establish estimations for the length of the second fundamental forms of hemi-slant submanifolds immersed as warped products in an lcK manifold.
	\section{Preliminaries}
	\allowdisplaybreaks
	\begin{defn}\label{def:lcK manifold}
		A Hermitian Manifold $(\widetilde{M}^{2n},J,g )$ is said to be a \textit{locally conformal Kähler} (l.c.K.) manifold if there exists an open cover $\{U_i\}_{i\in I}$ of $\widetilde{M}^{2n}$ and a family $\{f_i\}_{i\in I}$ of $C^\infty$ functions $f_i:U_i\to\mathbb{R}$ such that for each $i\in I$, the metric 
		\begin{equation}\label{eq:g_i defn}
			g_i=e^{-f_i}g|_{U_i} 
		\end{equation} 
		on $U_i$ is a Kähler metric.
	\end{defn}
	Given an l.c.K. manifold $(\widetilde{M}^{2n},J,g )$, let $U,V$ denote smooth sections of $T\widetilde{M}^{2n}$, then the local 1-forms $df_i$ glue up to a globally defined closed 1-form $\omega$, called the \textit{Lee form}, and it satisfies the following equation
	\begin{equation}\label{eq:d(Fund 2-form)}
		d\Omega=\Omega\wedge\omega
	\end{equation} 
	where $\Omega(U,V)=g (JU,V)$ is the fundamental 2-form associated to $(J,g )$.\par
	Denote by $\Theta$ the global closed 1-form defined as $\Theta=\omega\circ J$. Then, $\Theta$ is called the \textit{anti Lee form}.\par
	Denote by $B$ and $A$ the vector fields equivalent to $\omega$ and $\Theta$ respectively with respect to $g$, i.e. $\omega(U)=g(B,U)$ and $ \Theta(U)=g( A,U)$.\par
	$B$ and $A$ are respectively called the \textit{Lee vector field} and the \textit{anti Lee vector field}, and are related as 
	\begin{equation}\label{eq:Relation Lee-anti Lee vect}
		A=-JB
	\end{equation}
	Let $\overline{\nabla}$ denote the Levi-Civita connection of $(\widetilde{M}^{2n},g)$ and $\widetilde{\nabla}_i$ denote the Levi-Civita connection of the local metrics $g_i$ for all $i\in I$. Then $\widetilde{\nabla}_i$ glue up to a globally defined torsion-free linear connection $\widetilde{\nabla}$ on $\widetilde{M}^{2n}$ given by
	\begin{equation}\label{eq:Weyl Conn defn}
		\widetilde{\nabla}_UV=\overline{\nabla}_UV-\frac{1}{2}\left\lbrace\omega(U)V+\omega(V)U-g(U,V )B \right\rbrace  
	\end{equation}
	where $U,V\in T\widetilde{M}^{2n}$ and satisfying 
	\begin{equation}\label{eq:Weyl Conn and g}
		\widetilde{\nabla}g=\omega\otimes g
	\end{equation}
	$\widetilde{\nabla}$ is called the \textit{Weyl connection} of the l.c.K. manifold $(\widetilde{M}^{2n},J,g )$. As  $g_i$ are Kähler metrics, the almost complex structure $J$ is parallel with respect to the Weyl connection, i.e. $\widetilde{\nabla}J=0$. This gives 
	\begin{equation}\label{eq:Riemm Conn and J}
		\overline{\nabla}_UJV=J\overline{\nabla}_UV+\frac{1}{2}\left\lbrace \Theta(V)U-\omega(V)JU-g(U,V)A+\Omega(U,V)B\right\rbrace 
	\end{equation}
	Now as $\omega$ is a closed form on $\widetilde{M}^{2n}$, we have 
	\begin{equation}\label{eq:Derivative Lee Form Riemm Conn}
		(\overline{\nabla}_U\omega)V=(\overline{\nabla}_V\omega)U
	\end{equation}
	Hence using \eqref{eq:Riemm Conn and J} and \eqref{eq:Derivative Lee Form Riemm Conn} we have
	\begin{align*}
		(\overline{\nabla}_U\Theta)V		&=U\omega(JV)-\omega(\overline{\nabla}_UJV)+\frac{1}{2}\Theta(V)\omega(U)-\frac{1}{2}\omega(V)\Theta(U)+g(JU,V)||B||^2
	\end{align*}
	as $\omega(A)=g(B,A)=0$ from \eqref{eq:Relation Lee-anti Lee vect} and $\omega(B)=g(B,B)=||B||^2$\par\noindent
	Thus, we have 
	\begin{equation}\label{eq:Derivative Anti-Lee Form Riemm Conn}
		(\overline{\nabla}_U\Theta)V=(\overline{\nabla}_U\omega)V+\frac{1}{2}\Theta(V)\omega(U)-\frac{1}{2}\omega(V)\Theta(U)+g(JU,V)||B||^2
	\end{equation}
	Let $M^m$ be a Riemannian manifold isometrically immersed in an l.c.K. manifold $(\widetilde{M}^{2n},J,g )$. Let $U,V,W$ denote smooth sections of $TM^m$ and $\xi,\eta$ denote smooth sections of $T^\perp M^m$. \par
	The Gauss and Weingarten formulae with respect to the Riemannian connection of $\widetilde{M}^{2n}$ are given as 
	\begin{align}
		\overline{\nabla}_UV&=\nabla_UV+h(U,V)\label{eq:Gauss - Riemm Conn}\\
		\overline{\nabla}_U\xi&=-\mathfrak{A}_\xi U+\nabla^\perp_U \xi\label{eq:Weingarten - Riemm Conn}
	\end{align}
	where $h$ is the second fundamental form, $\mathfrak{A}$ is the shape operator and $\nabla, \nabla^\perp$ are respectively the induced connections in the tangent bundle and the normal bundle of $M^m$ with respect to $\overline{\nabla}$.\par  
	The Gauss and Weingarten formulae with respect to the Weyl connection of $\widetilde{M}^{2n}$ are given as 
	\begin{align}
		\widetilde{\nabla}_UV&=\hat{\nabla}_UV+\widetilde{h}(U,V)\label{eq:Gauss - Weyl Conn}\\
		\widetilde{\nabla}_U\xi&=-\widetilde{\mathfrak{A}}_\xi U+\widetilde{\nabla}^\perp_U \xi\label{eq:Weingarten - Weyl Conn}
	\end{align}
	where $\widetilde{h}$ is the second fundamental form, $\widetilde{\mathfrak{A}}$ is the shape operator and $\hat{\nabla}, \widetilde{\nabla}^\perp$ are respectively the induced connections in the tangent bundle and the normal bundle of $M^m$ with respect to $\widetilde{\nabla}$.\par
	Let $H$ denote the trace of $h$, then $H$ is called the mean curvature vector of $M^m$ in $(\widetilde{M}^{2n},J,g )$ and is a smooth section of $T^\perp M^m$. We say $M^m$ is a totally umbilic submanifold of $(\widetilde{M}^{2n},J,g )$, if $h(U,V)=g(U,V)H$. We say $M^m$ is a totally geodesic submanifold of $(\widetilde{M}^{2n},J,g )$, if $h(U,V)=0$.\par
	Let $B^T, B^N$ denote the tangential and normal components of the Lee vector field $B$ and let $A^T, A^N$ denote the tangential and normal components of the anti Lee vector field $A$.\par
	From \eqref{eq:Weyl Conn defn}, we have the following relations 
	\begin{align}
		\hat{\nabla}_UV&=\nabla_UV-\frac{1}{2}\left\lbrace\omega(U)V+\omega(V)U-g(U,V)B^T  \right\rbrace \label{eq:Relation induced conn M - Weyl and Riemm}\\
		\widetilde{h}(U,V)&=h(U,V)+\frac{1}{2}g(U,V)B^N \label{eq:Relation 2nd fund form - Weyl and Riemm}\\
		\widetilde{\mathfrak{A}}_\xi U&=\mathfrak{A}_\xi U+\frac{1}{2}\omega(\xi)U \label{eq:Relation shape operator - Weyl and Riemm}\\
		\widetilde{\nabla}^\perp_U\xi&=\nabla^\perp_U\xi-\frac{1}{2}\omega(U)\xi \label{eq:Relation induced normal conn - Weyl and Riemm}
	\end{align}
	Now define
	\begin{align}\label{eq:Definition P,F,t,f}
		JU&=PU+FU&J\xi=t\xi+f\xi
	\end{align}
	where $PU, t\xi$ and $FU, f\xi$ are respectively the tangential and normal parts. Then, we have
	\begin{equation}\label{eq:Identities P,F,t,f}
		\begin{aligned}
			P^2+tF&=-I&\hspace{1cm}f^2+Ft&=-I\\
			FP+fF&=0&\hspace{1cm}tf+Pt&=0
		\end{aligned}
	\end{equation}
	Now from \eqref{eq:Relation Lee-anti Lee vect} and \eqref{eq:Definition P,F,t,f} we have
	\begin{align}\label{eq:Relation A^T, A^N in terms of B^T, B^N}
		A^T&=-PB^T-tB^N&A^N&=-FB^T-fB^N
	\end{align}
	Define the covariant differentiation of $P$, $F$, $t$ and $f$ with respect to the Levi-Civita connection of $\widetilde{M}^{2n}$ as 
	\begin{equation}\label{eq:Definition Covar Diff Riem Conn P,F,t,f}
		\begin{aligned}
			(\overline{\nabla}_UP)V& =\nabla_UPV-P\nabla_UV\\
			(\overline{\nabla}_UF)V& =\nabla^\perp_UFV-F\nabla_UV\\
			(\overline{\nabla}_Ut)\xi& =\nabla_Ut\xi-t(\nabla^\perp_U\xi)\\
			(\overline{\nabla}_Uf)\xi& =\nabla^\perp_Uf\xi-f(\nabla^\perp_U\xi)
		\end{aligned}
	\end{equation}
	Similarly, define the covariant differentiation of $P$, $F$, $t$ and $f$ with respect to the Weyl connection of $\widetilde{M}^{2n}$ as 
	\begin{equation}\label{eq:Definition Covar Diff Weyl Conn P,F,t,f}
		\begin{aligned}
			(\widetilde{\nabla}_UP)V&=\hat{\nabla}_UPV-P\hat{\nabla}_UV\\
			(\widetilde{\nabla}_UF)V& =\widetilde{\nabla}^\perp_UFV-F\hat{\nabla}_UV\\
			(\widetilde{\nabla}_Ut)\xi&=\hat{\nabla}_Ut\xi-t\widetilde{\nabla}^\perp_U\xi\\
			(\widetilde{\nabla}_Uf)\xi&=\widetilde{\nabla}^\perp_Uf\xi-f\widetilde{\nabla}^\perp_U\xi
		\end{aligned}
	\end{equation}
	Then as $\widetilde{\nabla}J=0$, using \eqref{eq:Relation induced conn M - Weyl and Riemm}, \eqref{eq:Relation 2nd fund form - Weyl and Riemm}, \eqref{eq:Relation shape operator - Weyl and Riemm}, \eqref{eq:Relation induced normal conn - Weyl and Riemm} we have 
	\begin{equation}\label{eq:Relation Covar Diff Riem Conn P,F,t,f}
		\begin{aligned}
			(\overline{\nabla}_UP)V&=\mathfrak{A}_{FV}U+th(U,V)+\frac{1}{2}\left\lbrace \Theta(V)U-\omega(V)PU+g(PU,V)B^T-g(U,V)A^T\right\rbrace\\
			(\overline{\nabla}_UF)V&=fh(U,V)-h(U,PV)+\frac{1}{2}\left\lbrace g(PU,V)B^N-g(U,V)A^N-\omega(V)FU\right\rbrace\\
			(\overline{\nabla}_Ut)\xi&=\mathfrak{A}_{f\xi}U-P\mathfrak{A}_\xi U+\frac{1}{2}\left\lbrace g(FU,\xi)B^T-\omega(\xi)PU+\Theta(\xi)U \right\rbrace\\
			(\overline{\nabla}_Uf)\xi&=-h(U,t\xi)-F\mathfrak{A}_\xi U+\frac{1}{2}\left\lbrace g(FU,\xi)B^N-\omega(\xi)FU\right\rbrace 
		\end{aligned}
	\end{equation}
	Define the covariant derivative of the second fundamental form $h$ of the Riemannian connection $\overline{\nabla}$ as 
	\begin{equation}\label{eq:Definition Covar Diff Riem Conn h}
		(\overline{\nabla}_Uh)(V,W)=\nabla^\perp_Uh(V,W)-h(\nabla_UV,W)-h(V,\nabla_UW)
	\end{equation}
	Let $\overline{R}, R, R^\perp$ denote the curvature tensors associated to $\overline{\nabla}, \nabla, \nabla^\perp$ respectively. Then the Gauss, Codazzi and Ricci equations are respectively given by 
	\begin{align}
		g(\overline{R}(U,V)W,S)&=g(R(U,V)W,S)+g(h(V,S),h(U,W))-g(h(U,S),h(V,W))\label{eq:Gauss Eqn Riemm Conn}\\
		(\overline{R}(U,V)W)^\perp&=(\overline{\nabla}_Uh)(V,W)-(\overline{\nabla}_Vh)(U,W)\label{eq:Codazzi Eqn Riemm Conn}\\
		g(\overline{R}(U,V)\xi,\eta)&=g(R^\perp(U,V)\xi,\eta)-g([A_\xi,A_\eta]U,V)\label{eq: Ricci Eqn Riemm Conn}
	\end{align}
	Bishop and O'Neill \cite{WarpedProd-Def} defined warped product as
	\begin{defn}\label{defn:Warped Product Manifolds}
		Let $(M_1^{n_1},g_1)$ and $(M_2^{n_2},g_2)$ be Riemmanian manifolds and let \mbox{$\pi_1:M_1\times M_2\to M_1$} and \mbox{$\pi_2:M_1\times M_2\to M_2$} be the canonical projections. Let \mbox{$\lambda:M_1\to(0,\infty)$} be a smooth function. Then the warped product manifold \mbox{$(M,g)=M_1\times \,_\lambda M_2$} is defined as the manifold $M_1\times M_2$ equipped with the Riemannian metric
		\begin{equation}\label{eq:Warped Product Metric}
			g=\pi_1^\star g_1+\lambda^2\pi_2^\star g_2
		\end{equation}
	\end{defn}
	\noindent Warped product manifolds are a generalization of the usual product of two Riemannian manifolds. In fact we have the following characterisation theorem.
	\begin{theo}[\cite{Hiepko-WarpedProdCharacterisation}]\label{th:Characterisation - Warped Product}
		Let $(M^m,g)$ be a connected Riemannian manifold equipped with orthogonal, complementary, involutive distributions $\mathcal{D}_1$ and $\mathcal{D}_2$. Further let the leaves of $\mathcal{D}_1$ be totally geodesic and the leaves of $\mathcal{D}_2$ be extrinsic spheres in $M^m$, where by extrinsic spheres we mean totally umbilic submanifolds such that the mean curvature vector is parallel in the normal bundle. Then $(M^m,g)$ is locally a warped product \mbox{$(M,g)=M_1\times \,_\lambda M_2$}, where $M_1$ and $M_2$ respectively denote the leaves of $\mathcal{D}_1$ and $\mathcal{D}_2$ and \mbox{$\lambda:M_1\to(0,\infty)$} is a smooth function such that $\text{grad}(\ln\lambda)$ is the mean curvature vector of $M_2$ in $M$.\par
		Further, if $(M^m,g)$ is simply connected and complete, then $(M^m,g)$ is globally a warped product.
	\end{theo}
	\noindent For $(M_1^{n_1},g_1)$, $(M_2^{n_2},g_2)$ and $(M,g)$ denote respectively the Levi-Civita connections by $\nabla^1$, $\nabla^2$ and $\nabla$. Given any smooth function $\lambda:M_1\to\mathbb{R}$, let $\text{grad}(\lambda)$ denote the lift of the gradient vector field of $\lambda$ to $(M,g)$.
	\begin{theo}[\cite{Hiepko-WarpedProdCharacterisation}]\label{th:Expressions Connection Warped Product}
		Given a warped product manifold \mbox{$(M,g)= M_1\times \,_\lambda M_2$} of Riemmanian manifolds $(M_1^{n_1},g_1)$ and $(M_2^{n_2},g_2)$, we have for all $X,Y\in \mathcal{L}(M_1)$ and $Z,W\in \mathcal{L}(M_2)$,
		\begin{align}
			\nabla_XY&=\nabla^1_XY\label{eq:Connection Warped Product on M_1 x M_1}\\
			\nabla_XZ&=\nabla_ZX=X(\ln\lambda)Z\label{eq:Connection Warped Product on M_1 x M_2}\\
			\nabla_ZW&=\nabla^2_ZW-g(Z,W)\text{grad}(\ln\lambda)\label{eq:Connection Warped Product on M_2 x M_2}
		\end{align}
	\end{theo}
	\noindent It follows from Theorem \ref{th:Expressions Connection Warped Product} that $\mathcal{H}=-\text{grad}(\ln\lambda)$ is the mean curvature vector of $M_2$ in $M$.
	Let $M^m$ be a Riemannian manifold isometrically immersed in an l.c.K. manifold $(\widetilde{M}^{2n},J,g )$.\par
	$M^m$ is said to be a \textit{hemi-slant submanifold} if it admits two orthogonal complementary distributions $\mathcal{D}^\perp$ and $\mathcal{D}^\theta$, such that $\mathcal{D}^\perp$ is totally real, i.e. $J\mathcal{D}^\perp\subseteq T^\perp M^m$ and $\mathcal{D}^\theta$ is slant with slant angle $\theta\neq0,\frac{\pi}{2}$, i.e. $P^2Z=-\cos^2\theta Z$, for every smooth vector field $Z\in\mathcal{D}^\theta$.\par
	The tangent bundle and the normal bundle of a hemi-slant submanifold admits an orthogonal decomposition as
	\begin{align}\label{eq:Hemi Slant Tangent, Normal Bundle Ortho decomp }
		TM^m&=\mathcal{D}^\perp\oplus\mathcal{D}^\theta&T^\perp M^m&=J\mathcal{D}^\perp\oplus F\mathcal{D}^\theta\oplus\mu
	\end{align}
	where $\mu$ is the orthogonal complementary distribution of $J\mathcal{D}^\perp\oplus F\mathcal{D}^\theta$ in $T^\perp M^m$ and is an invariant subbundle of $T^\perp M^m$ with respect to $J$. It is easy to observe that,
	\begin{equation}\label{eq:Hemi-Slant - Image of D,D' under P,F,t}
		\begin{aligned}
			P\mathcal{D}^\perp&=\{0\}&\hspace{1cm}P\mathcal{D}^\theta&=\mathcal{D}^\theta&\hspace{1cm}F\mathcal{D}^\perp&=J\mathcal{D}^\perp\\
			t(J\mathcal{D}^\perp)&=\mathcal{D}^\perp&\hspace{1cm}t(F\mathcal{D}^\theta)&=\mathcal{D}^\theta&\hspace{1cm}t(\mu)&=\{0\}\\
			f(F\mathcal{D}^\perp)&=\{0\}&\hspace{1cm}f(F\mathcal{D}^\theta)&= F\mathcal{D}^\theta&\hspace{1cm}f(\mu)&=\mu
		\end{aligned}
	\end{equation}
	Let $M^m$ be a hemi-slant manifold isometrically immersed in an l.c.K. manifold $(\widetilde{M}^{2n},J,g )$ such that the totally real distribution $\mathcal{D}^\perp$ and the slant distribution $\mathcal{D}^\theta$ are both involutive. Let $M_\perp^{n_1}$ and $M_\theta^{2n_2}$ respectively denote the leaves of $\mathcal{D}^\perp$ and $\mathcal{D}^\theta$, where $n_1=\dim_\mathbb{R}\mathcal{D}^\perp$ and $2n_2=\dim_\mathbb{R}\mathcal{D}^\theta$. We say $M^m$ is a 
	\begin{itemize}
		\item \textit{mixed totally geodesic hemi-slant submanifold} if $h(\mathcal{D}^\perp,\mathcal{D}^\theta)=\{0\}$.
		\item \textit{hemi-slant product submanifold} if $M^m$ can be expressed locally as \mbox{$M_\perp\times M_\theta$}.
		\item \textit{hemi-slant warped product submanifold} if $M^m$ can be expressed locally as \mbox{$M_\perp\times\,_\lambda M_\theta$} for some smooth function $\lambda:M_\perp\to(0,\infty)$.
		\item \textit{warped product hemi-slant submanifold} if $M^m$ can be expressed locally as \mbox{$M_\theta\times\,_\lambda M_\perp$} for some smooth function $\lambda:M_\theta\to(0,\infty)$.
	\end{itemize}
	From here on we use $X,Y,X_1$ to denote smooth vector fields in $\mathcal{L}(M_\perp)$ and $Z,W$ to denote smooth vector fields in $\mathcal{L}(M_\theta)$.
	\begin{theo}\label{th:Hemi-Slant D_perp}\cite{Tastan_Sibel-HemiSlant-lcK}
		Let $M^m$ be a hemi-slant submanifold of an l.c.K. manifold $\widetilde{M}^{2n}$. Then 
		\begin{itemize}
			\item the totally real distribution $\mathcal{D}^\perp$ is involutive.
			\item the leaves of the totally real distribution $\mathcal{D}^\perp$ are totally geodesic in $M^m$ if and only if
			\begin{equation}\label{eq:hemi-slant D_perp totally geodesic}
				g(\mathfrak{A}_{JX}Z-\mathfrak{A}_{FZ}X,Y)=\frac{1}{2}\omega(PZ)g(X,Y)
			\end{equation}
			\item the leaves of the totally real distribution $\mathcal{D}^\perp$ are totally umbilic in $M^m$ if and only if 
			\begin{equation}\label{eq:hemi-slant D_perp totally umbilic}
				g(\mathfrak{A}_{JX}Z-\mathfrak{A}_{FZ}X,Y)=\left(\frac{1}{2}\omega(PZ)+g(\mathcal{H},PZ)\right)g(X,Y)
			\end{equation}
			for some smooth vector field $\mathcal{H}\in\mathcal{D}^\theta$.
		\end{itemize}
	\end{theo}
	\begin{theo}\label{th:Hemi-Slant D_theta}\cite{Tastan_Sibel-HemiSlant-lcK}
		Let $M^m$ be a hemi-slant submanifold of an l.c.K. manifold $\widetilde{M}^{2n}$. Then 
		\begin{itemize}
			\item the slant distribution $\mathcal{D}^\theta$ is involutive if and only if 
			\begin{equation}\label{eq:hemi-slant D_theta involutive}
				g(\mathfrak{A}_{FPZ}X,W)+g(\nabla^\perp_WFZ,JX)=g(\mathfrak{A}_{FPW}X,Z)+g(\nabla^\perp_ZFW,JX)
			\end{equation}
			\item the leaves of the slant distribution $\mathcal{D}^\theta$ are totally geodesic in $M^m$ if and only if 
			\begin{equation}\label{eq:hemi-slant D_theta totally geodesic}
				\omega(\mathcal{D}^\perp)=\{0\}\text{ and }g(\mathfrak{A}_{FPZ}X,W)+g(\nabla^\perp_WFZ,JX)=0
			\end{equation}
			\item the leaves of the slant distribution $\mathcal{D}^\theta$ are totally umbilic in $M^m$ if and only if 
			\begin{equation}\label{eq:hemi-slant D_theta totally umbilic}
				g(\mathfrak{A}_{FPZ}X,W)+g(\nabla^\perp_WFZ,JX)=\sin^2\theta\left(\frac{1}{2}\omega(X)+g(\mathcal{H},X)\right)g(Z,W)
			\end{equation}
			for some smooth vector field $\mathcal{H}\in\mathcal{D}^\perp$.
		\end{itemize}
	\end{theo}
	\noindent \textbf{Notations:} Let $\mathcal{D}^\perp$ and $\mathcal{D}^\theta$ be the totally real and slant distributions on a hemi-slant submanifold $M^m$ of an lcK manifold $\widetilde{M}^{2n}$ such that both distributions are involutive and let $M_\perp$ and $M_\theta$ respectively denote the leaves of the distributions $\mathcal{D}^\perp$ and $\mathcal{D}^\theta$ respectively. Then $\mathcal{D}^\perp(p,q)=T_{(p,q)}(M_\perp\times\{q\})$ and $\mathcal{D}^\theta(p,q)=T_{(p,q)}(\{p\}\times M_\theta)$. Let $\mathcal{L}(M_\perp)$ and $\mathcal{L}(M_\theta)$ respectively denote the set of lifts of vector fields from $M_\perp$ and $M_\theta$ to $M$. Then $X\in\mathcal{L}(M_\perp)$ if and only if $X|_{\{p\}\times M_\theta}$ is constant for every $p\in M_\perp$. Similarly, $Z\in\mathcal{L}(M_\theta)$ if and only if $Z|_{M_\perp\times\{q\}}$ is constant for every $q\in M_\theta$. Also, if \mbox{$\pi_\perp:M_\perp\times M_\theta\to M_\perp$} and \mbox{$\pi_\theta:M_\perp\times M_\theta\to M_\theta$} are the canonical projections, we have $d\pi_\perp(\mathcal{L}(M_\perp))=TM_\perp$ and $d\pi_\theta(\mathcal{L}(M_\theta))=TM_\theta$. It is clear that a general vector field in $\mathcal{D}^\perp$ (respectively $\mathcal{D}^\theta$) need not be in $\mathcal{L}(M_\perp)$ (respectively $\mathcal{L}(M_\theta)$).
	
	\section{Hemi-Slant Warped Product Submanifolds of l.c.K. manifolds}
	\begin{lem}\label{lem:Identities For Hemi-Slant Warped Product Submanifolds}
		Given a hemi-slant warped product submanifold \mbox{$M= M_\perp\times \,_{\lambda}M_\theta$} in an lcK manifold $(\widetilde{M}^{2n},J,g )$, we have for all $X,Y,X_1\in \mathcal{L}(M_\perp)$ and $Z,W\in \mathcal{L}(M_\theta)$, 
		\begin{align}
			g(h(X,Z),JY)&=g(h(Y,Z),JX)\label{eq:g(h(D_perp,D_theta),JD_perp) Symmetric Hemi-Slant Warped Product Submanifolds}\\
			g(h(X,Z),FW)&=g(h(X,W),FZ)\label{eq:g(h(D_perp,D_theta),FD_theta) Symmetric Hemi-Slant Warped Product Submanifolds}\\
			g(h(Z,W),JX)&=g(h(X,Z),FW)+\frac{1}{2}g(Z,W)g(JB,X)\label{eq:Identity g(h(D_theta,D_theta),JD_perp) and g(h(D_perp,D_theta),FD_theta) Hemi-Slant Warped Product Submanifolds}\\
			X(\ln\lambda)&=\frac{1}{2}g(B,X)\label{eq:Identity g(B,D_perp) Hemi-Slant Warped Product Submanifolds}\\
			g(h(X,Y),JX_1)&=g(h(X,X_1),JY)-\frac{1}{2}g(X,Y)g(B,JX_1)+\frac{1}{2}g(X,X_1)g(B,JY)\label{eq:Identity g(h(D_perp,D_perp),FD_perp) Hemi-Slant Warped Product Submanifolds}
		\end{align}
	\end{lem}
	\begin{proof}
		For all $X,Y,X_1\in \mathcal{L}(M_\perp)$ and $Z,W\in \mathcal{L}(M_\theta)$, we have using \eqref{eq:Riemm Conn and J} and \eqref{eq:Connection Warped Product on M_1 x M_2},
		\begin{align*}
			g(h(X,Z),JY)&=g(\overline{\nabla}_ZX,JY)=-g(J\overline{\nabla}_ZX,Y)=-g(\overline{\nabla}_ZJX,Y)=g(\mathfrak{A}_{JX}Z,Y)\\
			&=g(h(Y,Z),JX)
			\intertext{which implies \eqref{eq:g(h(D_perp,D_theta),JD_perp) Symmetric Hemi-Slant Warped Product Submanifolds}. Similarly,}
			g(h(X,Z),FW)&=g(\overline{\nabla}_XZ,JW)-g(\overline{\nabla}_XZ,PW)\\
			&=-g(J\overline{\nabla}_XZ,W)-g(\nabla_XZ,PW)\\
			&=-g(\overline{\nabla}_XJZ,W)-X(\ln\lambda)g(Z,PW)\\
			&=-X(\ln\lambda)g(PZ,W)+g(\mathfrak{A}_{FZ}X,W)-X(\ln\lambda)g(Z,PW)
			\intertext{which implies \eqref{eq:g(h(D_perp,D_theta),FD_theta) Symmetric Hemi-Slant Warped Product Submanifolds}. Repeating the above calculation, we have}
			g(h(X,Z),FW)&=-g(J\overline{\nabla}_ZX,W)-g(\nabla_ZX,PW)\\
			&=-g(\overline{\nabla}_ZJX,W)-\frac{1}{2}g(JB,X)g(Z,W)-\frac{1}{2}g(B,X)g(JZ,W)-X(\ln\lambda)g(Z,PW)\\
			&=g(\mathfrak{A}_{JX}Z,W)-\frac{1}{2}g(JB,X)g(Z,W)-\frac{1}{2}g(B,X)g(PZ,W)-X(\ln\lambda)g(Z,PW)
			\intertext{Using \eqref{eq:g(h(D_perp,D_theta),FD_theta) Symmetric Hemi-Slant Warped Product Submanifolds} and comparing symmetric and skew symmetric terms in $Z$ and $W$ we have,}
			g(h(X,Z),FW)&=g(h(Z,W),JX)-\frac{1}{2}g(JB,X)g(Z,W)
			\intertext{which proves \eqref{eq:Identity g(h(D_theta,D_theta),JD_perp) and g(h(D_perp,D_theta),FD_theta) Hemi-Slant Warped Product Submanifolds} and}
			0&=\left( X(\ln\lambda)-\frac{1}{2}g(B,X)\right) g(PZ,W)
			\intertext{which proves \eqref{eq:Identity g(B,D_perp) Hemi-Slant Warped Product Submanifolds}. Finally,}
			g(h(X,Y),JX_1)&=g(\overline{\nabla}_XY,JX_1)\\
			&=-g(J\overline{\nabla}_XY,X_1)\\
			&=-g(\overline{\nabla}_XJY,X_1)-\frac{1}{2}g(X,X_1)g(JB,Y)+\frac{1}{2}g(X,Y)g(JB,X_1)\\
			&=g(\mathfrak{A}_{JY}X,X_1)-\frac{1}{2}g(X,X_1)g(JB,Y)+\frac{1}{2}g(X,Y)g(JB,X_1)
		\end{align*}
		which gives \eqref{eq:Identity g(h(D_perp,D_perp),FD_perp) Hemi-Slant Warped Product Submanifolds}.	
	\end{proof}
	\begin{rem}\label{rem:Local Orthonormal Basis of Hemi-Slant Warped Product Submanifolds}
		Given a hemi-slant warped product submanifold $M_\perp\times\,_\lambda M_\theta$ of an l.c.K manifold $\widetilde{M}^{2n}$, let $\{X_i\}_{i=1}^{p}$ and $\{Z_j,\beta PZ_j\}_{j=1}^{q}$ respectively be local orthonormal frames of $TM_\perp$ and $TM_\theta$. Then a local orthonormal frame of $\widetilde{M}^{2n}$ is
		\[\begin{aligned}
			\left\lbrace \widehat{X_i}=X_i\right\rbrace\cup\left\lbrace\widehat{Z_j}=\frac{Z_j}{\lambda},\widehat{PZ_j}=\frac{\beta PZ_j}{\lambda}\right\rbrace\cup\left\lbrace \widehat{JX_i}=JX_i\right\rbrace\\\cup\left\lbrace\widehat{FZ_j}=\frac{\alpha FZ_j}{\lambda},\widehat{FPZ_j}=\frac{\alpha\beta FPZ_j}{\lambda}\right\rbrace\cup\left\lbrace\widehat{\xi_k},\widehat{J\xi_k}\right\rbrace
		\end{aligned}\]
		where $\alpha=\csc\theta$, $\beta=\sec\theta$ and
		\begin{align*}
			&\left\{\widehat{X_i}:1\leq i\leq n_1\right\}\text{ is an orthonormal basis of }\mathcal{D}^\perp\\
			&\left\{\widehat{Z_j}, \widehat{PZ_j}:1\leq j\leq n_2\right\}\text{ is an orthonormal basis of }\mathcal{D}^\theta\\
			&\left\{\widehat{JX_i}:1\leq i\leq n_1\right\}\text{ is an orthonormal basis of }J\mathcal{D}^\perp\\
			&\left\{\widehat{FZ_j}, \widehat{FPZ_j}:1\leq j\leq n_2\right\}\text{ is an orthonormal basis of }F\mathcal{D}^\theta\\
			&\left\{\widehat{\xi_k},\widehat{J\xi_k}:1\leq k\leq \frac{n-n_1-2n_2}{2}\right\}\text{ is an orthonormal basis of }\mu
		\end{align*}
		However, while $Z_j,\beta PZ_j\in\mathcal{L}(M_\theta)$ we have $\widehat{Z_j},\widehat{PZ_j}\notin\mathcal{L}(M_\theta)$ in general, as $\lambda$ is a function on $M_\perp$. Also, note that
		\begin{align*}
			J\left(\widehat{Z_j}\right)&=J\left(\frac{Z_j}{\lambda}\right)=\frac{PZ_j}{\lambda}+\frac{FZ_j}{\lambda}=\cos\theta\widehat{PZ_j}+\sin\theta\widehat{FZ_j}\\
			J\left(\widehat{PZ_j}\right)&=J\left(\sec\theta\frac{PZ_j}{\lambda}\right)=\frac{\sec\theta P^2Z_j}{\lambda}+\frac{\sec\theta FPZ_j}{\lambda}=-\cos\theta\widehat{Z_j}+\sin\theta\widehat{FPZ_j}
		\end{align*}
	\end{rem}
	\noindent We now give a characterisation for hemi-slant warped product submanifolds of l.c.K. manifolds.
	\begin{theo}\label{th:Hemi-Slant Warped Prod Characterisation - M_perp x_lambda M_theta}
		Let $M^m$ be a hemi-slant submanifold of an l.c.K. manifold $\widetilde{M}^{2n}$. Then the following are equivalent 
		\begin{enumerate}
			\item $M^m$ is a hemi-slant warped product submanifold \mbox{$M_\perp\times\,_\lambda M_\theta$} of $\widetilde{M}^{2n}$\label{Hemi-Slant Warped Prod Characterisation}
			\item For every $X,Y\in\mathcal{L}(M_\perp)$ and $Z,W\in\mathcal{L}(M_\theta)$ we have
			\begin{equation}\label{eq:Hemi-Slant Warped Prod Characterisation 1}
				\begin{aligned}
					g(\mathfrak{A}_{JX}Z-\mathfrak{A}_{FZ}X,Y)&=\frac{1}{2}\omega(PZ)g(X,Y)\\
					g(\mathfrak{A}_{FPZ}X,W)+g(\nabla^\perp_WFZ,JX)&=\sin^2\theta\left(\frac{1}{2}\omega(X)-X(\ln\lambda)\right)g(Z,W)
				\end{aligned}
			\end{equation}
			for some smooth function $\lambda:M_\perp\to(0,\infty)$.\label{Hemi-Slant Warped Prod Characterisation - Condition 1}
			\item For every $X\in\mathcal{L}(M_\perp)$ and $Z\in\mathcal{L}(M_\theta)$ we have
			\begin{equation}\label{eq:Hemi-Slant Warped Prod Characterisation 2}
				\nabla_XZ=\nabla_ZX=\frac{1}{2}\omega(X)Z
			\end{equation}\label{Hemi-Slant Warped Prod Characterisation - Condition 2}
		\end{enumerate}
		Also, in this case we have the mean curvature vector $\mathcal{H}$ of $M_\theta$ in $M^m$ is
		\begin{equation}\label{eq:Warping Function Hemi-Slant Warped Prod Characterisation}
			\mathcal{H}=-\text{grad}(\ln\lambda)=-\frac{1}{2}B|_{\mathcal{D}^\perp}
		\end{equation}
		where $B|_{\mathcal{D}^\perp}$ is the component of $B$ along $\mathcal{D}^\perp$.
	\end{theo}
	\begin{proof}
		\mbox{\eqref{Hemi-Slant Warped Prod Characterisation}$\Leftrightarrow$\eqref{Hemi-Slant Warped Prod Characterisation - Condition 1}} This follows from Theorem \ref{th:Hemi-Slant D_perp}, Theorem \ref{th:Hemi-Slant D_theta} and the fact that $\nabla\ln\lambda\in\mathcal{L}(M_\perp)$ which implies for all $X\in \mathcal{L}(M_\perp)$ and $Z\in \mathcal{L}(M_\theta)$
		\begin{align*}
			g(\nabla_Z(\nabla\ln\lambda),X)&=ZX(\ln\lambda)-g(\nabla\ln\lambda,\nabla_ZX)\\
			&=[Z,X](\ln\lambda)-\nabla_ZX(\ln\lambda)\hspace{0.4cm}(\text{as }Z(\ln\lambda)=0)\\
			&=-\nabla_XZ(\ln\lambda)\\
			&=g(Z,\nabla_X(\nabla\ln\lambda))\\
			&=0\hspace{0.4cm}
		\end{align*}
		as $\mathcal{D}^\perp$ is totally geodesic.	Also, \eqref{eq:Warping Function Hemi-Slant Warped Prod Characterisation} follows from Lemma \ref{lem:Identities For Hemi-Slant Warped Product Submanifolds} \eqref{eq:Identity g(B,D_perp) Hemi-Slant Warped Product Submanifolds}.\par\noindent
		\mbox{\eqref{Hemi-Slant Warped Prod Characterisation}$\Leftrightarrow$\eqref{Hemi-Slant Warped Prod Characterisation - Condition 2}} Let \mbox{$M=M_\perp\times\,_\lambda M_\theta$} be a hemi-slant warped product submanifold. Then \eqref{eq:Hemi-Slant Warped Prod Characterisation 2} and \eqref{eq:Warping Function Hemi-Slant Warped Prod Characterisation} follow from \eqref{eq:Connection Warped Product on M_1 x M_2} and Lemma \ref{lem:Identities For Hemi-Slant Warped Product Submanifolds} \eqref{eq:Identity g(B,D_perp) Hemi-Slant Warped Product Submanifolds}.\par
		Conversely, let $M^m$ be a hemi-slant submanifold of an l.c.K. manifold $\widetilde{M}^{2n}$ such that \eqref{eq:Hemi-Slant Warped Prod Characterisation 2} holds. Then for all $X,Y\in \mathcal{L}(M_\perp)$ and $Z,W\in \mathcal{L}(M_\theta)$ we have
		\begin{align*}
			g([X,Y],Z)&=g(\nabla_XY-\nabla_YX,Z)\\
			&=-g(\nabla_XZ,Y)+g(\nabla_YZ,X)\\
			&=0
			\intertext{which implies $\mathcal{D}^\perp$ is involutive.}
			g(\nabla_XY,Z)&=-g(\nabla_XZ,Y)=0
			\intertext{which implies leaves of $\mathcal{D}^\perp$ are totally geodesic in $M$.}
			g([Z,W],X)&=g(\nabla_ZW-\nabla_WZ,X)\\
			&=-g(\nabla_ZX,W)+g(\nabla_WX,Z)\\
			&=-\frac{1}{2}\omega(X)g(Z,W)+\frac{1}{2}\omega(X)g(W,Z)=0
			\intertext{which implies $\mathcal{D}^\theta$ is involutive.}
			g(\nabla_ZW,X)&=-g(\nabla_ZX,W)\\
			&=-\frac{1}{2}\omega(X)g(Z,W)\\
			&=-\frac{1}{2}g(Z,W)g(B^T,X)
			\intertext{which implies leaves of $\mathcal{D}^\theta$ are totally umbilical in $M$ with mean curvature vector $-\frac{1}{2}B|_{\mathcal{D}^\perp}$.}
			g\left(\nabla_ZB|_{\mathcal{D}^\perp},X\right)&=\frac{1}{2}\omega\left(B|_{\mathcal{D}^\perp}\right)g(Z,X)=0
		\end{align*}
		which implies $B|_{\mathcal{D}^\perp}$ is parallel in the normal bundle of $M_\theta$ in $M$.\par
		Hence by Theorem \ref{th:Characterisation - Warped Product} we have \mbox{$M=M_\perp\times\,_\lambda M_\theta$} is a hemi-slant warped product submanifold.
	\end{proof}
	We conclude our study of hemi-slant warped product submanifolds of l.c.K. manifolds by giving an inequality for the norm of the second fundamental form.
	\begin{theo}\label{th:Length h - Hemi-Slant Warped Product}
		Let \mbox{$M= M_\perp\times \,_\lambda M_\theta$} be a hemi-slant warped product submanifold in an lcK manifold $(\widetilde{M}^{2n},J,g )$. Then the norm of the second fundamental form satisfies the inequality
		\begin{align}\label{eq:Length h - Hemi-Slant Warped Product}
			||h||^2\geq\:&
			\frac{(n_1+n_2-1)}{2}\|B|_{J\mathcal{D}^\perp}\|^2+2g\left(H_{\mathcal{D}^\perp}|_{J\mathcal{D}^\perp},B|_{J\mathcal{D}^\perp}\right)+g\left(H_{\mathcal{D}^\theta}|_{J\mathcal{D}^\perp},B|_{J\mathcal{D}^\perp}\right)-K
		\end{align}
		where $n_1=\dim_{\mathbb{R}}\mathcal{D}^\perp$, $2n_2=\dim_{\mathbb{R}}\mathcal{D}^\theta$, $H_{\mathcal{D}^\perp}$ and $H_{\mathcal{D}^\theta}$ are respectively the components of the mean curvature vector $H$ of $M$ in $\widetilde{M}^{2n}$ along $\mathcal{D}^\perp$ and $\mathcal{D}^\theta$ and given any orthonormal basis $\left\{\widehat{X_i}:1\leq i\leq n_1\right\}$ of $\mathcal{D}^\perp$, \mbox{$K=2\sum_{i}g\left( h\left(\widehat{X_i},\widehat{X_i}\right),\widehat{JX_i}\right)g\left(B,\widehat{JX_i}\right)$}  .\par
		\noindent If equality holds then we have
		\begin{itemize}
			\item $\text{Image}(h)\subseteq(J\mathcal{D}^\perp\oplus F\mathcal{D}^\theta)$, and
			\item $M_\theta$ is totally umbilical in $\widetilde{M}^{2n}$ (with mean curvature vector $\mathcal{H}=-\frac{1}{2}B|_{\mathcal{D}^\perp}$) if and only if, $M$ is mixed-totally geodesic in $\widetilde{M}^{2n}$.
		\end{itemize}
	\end{theo}
	\begin{proof}
		\begin{align*}
			||h||^2=&\left\|h(\mathcal{D}^\perp,\mathcal{D}^\perp)\big|_{J\mathcal{D}^\perp}\right\|^2+\left\|h(\mathcal{D}^\perp,\mathcal{D}^\theta)\big|_{J\mathcal{D}^\perp}\right\|^2+\left\|h(\mathcal{D}^\theta,\mathcal{D}^\theta)\big|_{J\mathcal{D}^\perp}\right\|^2+\left\|h(\mathcal{D}^\perp,\mathcal{D}^\perp)\big|_{F\mathcal{D}^\theta}\right\|^2+\left\|h(\mathcal{D}^\perp,\mathcal{D}^\theta)\big|_{F\mathcal{D}^\theta}\right\|^2\\
			&+\left\|h(\mathcal{D}^\theta,\mathcal{D}^\theta)\big|_{F\mathcal{D}^\theta}\right\|^2+\left\|h(\mathcal{D}^\perp,\mathcal{D}^\perp)\big|_{\mu}\right\|^2+\left\|h(\mathcal{D}^\perp,\mathcal{D}^\theta)\big|_{\mu}\right\|^2+\left\|h(\mathcal{D}^\theta,\mathcal{D}^\theta)\big|_{\mu}\right\|^2
		\end{align*}
		From \eqref{eq:Identity g(h(D_perp,D_perp),FD_perp) Hemi-Slant Warped Product Submanifolds} and Remark \ref{rem:Local Orthonormal Basis of Hemi-Slant Warped Product Submanifolds} we have
		\begin{align*}
			\left\|h(\mathcal{D}^\perp,\mathcal{D}^\perp)\big|_{J\mathcal{D}^\perp}\right\|^2
			=&\sum_{i}g\left( h\left(\widehat{X_i},\widehat{X_i}\right),\widehat{JX_i}\right)^2
			+\sum_{i\ne j}\left\lbrace
			g\left( h\left(\widehat{X_i},\widehat{X_i}\right),\widehat{JX_j}\right)^2
			+2g\left( h\left(\widehat{X_i},\widehat{X_j}\right),\widehat{JX_i}\right)^2
			\right\rbrace\\
			&+\sum_{i\ne j \ne k}g\left( h\left(\widehat{X_i},\widehat{X_j}\right),\widehat{JX_k}\right)^2\\
			\geq\:&\frac{2}{\lambda^6}\sum_{i\ne j}g\left( h\left(X_i,X_j\right),JX_i\right)^2\\
			=&\frac{2}{\lambda^6}\sum_{i\ne j}
			\biggl\{g\left( h\left(X_i,Z_i\right),JX_j\right)
			+\frac{1}{2}g\left(X_i,X_i\right)g\left(B,JX_j\right)\biggr\}^2\\
			=&2\sum_{i\ne j}g\left( h\left(\widehat{X_i},\widehat{X_i}\right),\widehat{JX_j}\right)^2
			+\frac{1}{2}\sum_{i\ne j}g\left(\widehat{X_i},\widehat{X_i}\right)^2g\left(B,\widehat{JX_j}\right)^2\\
			&+2\sum_{i\ne j}g\left( h\left(\widehat{X_i},\widehat{X_i}\right),\widehat{JX_j}\right)g\left(\widehat{X_i},\widehat{X_i}\right)g\left(B,\widehat{JX_j}\right)\\
			\geq\:&\frac{p-1}{2}\biggl\|B\big|_{J\mathcal{D}^\perp}\biggr\|^2
			+2g\left(\sum_{i} \left. h\left(\widehat{X_i},\widehat{X_i}\right)\right|_{J\mathcal{D}^\perp},B\big|_{J\mathcal{D}^\perp}\right)\\
			&-2\sum_{i}g\left( h\left(\widehat{X_i},\widehat{X_i}\right),\widehat{JX_i}\right)g\left(B,\widehat{JX_i}\right)\\
			=&\frac{p-1}{2}\biggl\|B\big|_{J\mathcal{D}^\perp}\biggr\|^2
			+2g\left(H_{\mathcal{D}^\perp}|_{J\mathcal{D}^\perp},B|_{J\mathcal{D}^\perp}\right)-K
		\end{align*}
		From \eqref{eq:Identity g(h(D_theta,D_theta),JD_perp) and g(h(D_perp,D_theta),FD_theta) Hemi-Slant Warped Product Submanifolds} and Remark \ref{rem:Local Orthonormal Basis of Hemi-Slant Warped Product Submanifolds} we have
		\begin{align*}
			g\left( h\left(\widehat{Z_p},\widehat{Z_q}\right),\widehat{JX_i}\right)
			=&\frac{1}{\lambda^2}g\left( h\left(Z_p,Z_q\right),JX_i\right)
			=\frac{1}{\lambda^2}\left\lbrace g\left( h\left(X_i,Z_p\right),FZ_q\right)
			+\frac{1}{2}\lambda^2\delta_{pq}g\left( JB,X_i\right) \right\rbrace\\
			=&\sin\theta g\left( h\left(\widehat{X_i},\widehat{Z_p}\right),\widehat{FZ_q}\right)
			+\frac{1}{2}\delta_{pq} g\left( JB,\widehat{X_i}\right)\\
			g\left( h\left(\widehat{Z_p},\widehat{PZ_q}\right),\widehat{JX_i}\right)
			=&\frac{\sec\theta}{\lambda^2}g\left( h\left(Z_p,PZ_q\right),JX_i\right)
			=\frac{\sec\theta}{\lambda^2} g\left( h\left(X_i,Z_p\right),FPZ_q\right)\\
			=&\sin\theta g\left( h\left(\widehat{X_i},\widehat{Z_p}\right),\widehat{FPZ_q}\right)\\
			g\left( h\left(\widehat{PZ_p},\widehat{Z_q}\right),\widehat{JX_i}\right)
			=&\sin\theta g\left( h\left(\widehat{X_i},\widehat{PZ_p}\right),\widehat{FZ_q}\right)\\
			g\left( h\left(\widehat{PZ_p},\widehat{PZ_q}\right),\widehat{JX_i}\right)
			=&\frac{\sec^2\theta}{\lambda^2}g\left( h\left(PZ_p,PZ_q\right),JX_i\right)\\
			=&\frac{\sec^2\theta}{\lambda^2}\left\lbrace g\left( h\left(X_i,PZ_p\right),FPZ_q\right)
			+\frac{1}{2}\lambda^2\cos^2\theta\delta_{pq}g\left( JB,X_i\right) \right\rbrace\\
			=&\sin\theta g\left( h\left(\widehat{X_i},\widehat{PZ_p}\right),\widehat{FPZ_q}\right)
			+\frac{1}{2}\delta_{pq} g\left( JB,\widehat{X_i}\right)
		\end{align*}
		which implies
		\begin{align*}
			\left\|h(\mathcal{D}^\theta,\mathcal{D}^\theta)\big|_{J\mathcal{D}^\perp}\right\|^2
			=&\sum_{i,p,q}\left\lbrace
			g\left( h\left(\widehat{Z_p},\widehat{Z_q}\right),\widehat{JX_i}\right)^2
			+g\left( h\left(\widehat{Z_p},\widehat{PZ_q}\right),\widehat{JX_i}\right)^2
			+g\left( h\left(\widehat{PZ_p},\widehat{Z_q}\right),\widehat{JX_i}\right)^2 \right.\\
			&\left.+g\left( h\left(\widehat{PZ_p},\widehat{PZ_q}\right),\widehat{JX_i}\right)^2
			\right\rbrace\\
			=&\sin^2\theta\sum_{i,p,q}\left\lbrace
			g\left( h\left(\widehat{X_i},\widehat{Z_p}\right),\widehat{FZ_q}\right)^2
			+g\left( h\left(\widehat{X_i},\widehat{Z_p}\right),\widehat{FPZ_q}\right)^2
			+g\left( h\left(\widehat{X_i},\widehat{PZ_p}\right),\widehat{FZ_q}\right)^2 \right.\\
			&\left.+g\left( h\left(\widehat{X_i},\widehat{PZ_p}\right),\widehat{FPZ_q}\right)^2
			\right\rbrace
			-\frac{1}{2}\sum_{i,p}g\left( JB,\widehat{X_i}\right)^2
			+\sum_{i,p}\left\lbrace 
			g\left( h\left(\widehat{Z_p},\widehat{Z_p}\right),\widehat{JX_i}\right)g\left( JB,\widehat{X_i}\right)\right.\\
			&\left.+g\left( h\left(\widehat{PZ_p},\widehat{PZ_p}\right),\widehat{JX_i}\right)g\left( JB,\widehat{X_i}\right)
			\right\rbrace\\
			=&\sin^2\theta\left\|h(\mathcal{D}^\perp,\mathcal{D}^\theta)\big|_{F\mathcal{D}^\theta}\right\|^2
			-g\left(\left.\sum_p \left\lbrace h\left(\widehat{Z_p},\widehat{Z_p}\right)+ h\left(\widehat{PZ_p},\widehat{PZ_p}\right)\right\rbrace\right|_{J\mathcal{D}^\perp},B|_{J\mathcal{D}^\perp}\right)\\
			&-\frac{2q}{4}\biggl\|B\big|_{J\mathcal{D}^\perp}\biggr\|^2\\
			=&\sin^2\theta\left\|h(\mathcal{D}^\perp,\mathcal{D}^\theta)\big|_{F\mathcal{D}^\theta}\right\|^2
			-g\left(H_{\mathcal{D}^\theta}|_{J\mathcal{D}^\perp},B|_{J\mathcal{D}^\perp}\right)
			-\frac{q}{2}\biggl\|B\big|_{J\mathcal{D}^\perp}\biggr\|^2
		\end{align*}
		Combining we have \eqref{eq:Length h - Hemi-Slant Warped Product}.\par
		If equality holds in \eqref{eq:Length h - Hemi-Slant Warped Product}, then the only non-zero components of $||h||$ are  $\left\|h(\mathcal{D}^\perp,\mathcal{D}^\perp)|_{J\mathcal{D}^\perp}\right\|^2$, $\left\|h(\mathcal{D}^\perp,\mathcal{D}^\theta)|_{F\mathcal{D}^\theta}\right\|^2$ and $\left\|h(\mathcal{D}^\theta,\mathcal{D}^\theta)|_{J\mathcal{D}^\perp}\right\|^2$. Also, from the above computations we have, $\left\|h(\mathcal{D}^\perp,\mathcal{D}^\theta)|_{F\mathcal{D}^\theta}\right\|^2=0$ if and only if $\left\|h(\mathcal{D}^\theta,\mathcal{D}^\theta)|_{J\mathcal{D}^\perp}\right\|^2=0$. Hence, the result follows.
	\end{proof}
	
	\section{Warped Product Hemi-Slant Submanifolds of l.c.K. manifolds}
	\begin{lem}\label{lem:Identities For Warped Product Hemi-Slant Submanifolds}
		Given a warped product hemi-slant submanifold \mbox{$M= M_\theta\times \,_{\lambda}M_\perp$} in an lcK manifold $(\widetilde{M}^{2n},J,g )$, we have for all $X,Y\in \mathcal{L}(M_\perp)$ and $Z,W\in \mathcal{L}(M_\theta)$, 
		\begin{align}
			g(h(X,Z),JY)&=g(h(Y,Z),JX)\label{eq:g(h(D_perp,D_theta),JD_perp) Symmetric Warped Product Hemi-Slant Submanifolds}\\
			g(h(X,Z),FW)&=g(h(X,W),FZ)\label{eq:g(h(D_perp,D_theta),FD_theta) Symmetric Warped Product Hemi-Slant Submanifolds}\\
			g(h(Z,W),JX)&=g(h(X,Z),FW)+\frac{1}{2}g(Z,W)g(JB,X)\label{eq:Identity g(h(D_theta,D_theta),JD_perp) and g(h(D_perp,D_theta),FD_theta) Warped Product Hemi-Slant Submanifolds}\\
			g(B,X)&=0\label{eq:Identity g(B,D_perp) Warped Product Hemi-Slant Submanifolds}\\
			g(h(X,Y),JX_1)&=g(h(X,X_1),JY)-\frac{1}{2}g(X,Y)g(B,JX_1)+\frac{1}{2}g(X,X_1)g(B,JY)\label{eq:Identity g(h(D_perp,D_perp),FD_perp) Warped Product Hemi-Slant Submanifolds}
		\end{align}
	\end{lem}
	\begin{proof}
		For all $X,Y\in \mathcal{L}(M_\perp)$ and $Z,W\in \mathcal{L}(M_\theta)$, we have using \eqref{eq:Riemm Conn and J} and \eqref{eq:Connection Warped Product on M_1 x M_2},
		\begin{align*}
			g(h(X,Z),JY)&=g(\overline{\nabla}_ZX,JY)=-g(J\overline{\nabla}_ZX,Y)=-g(\overline{\nabla}_ZJX,Y)=g(\mathfrak{A}_{JX}Z,Y)\\
			&=g(h(Y,Z),JX)
			\intertext{which implies \eqref{eq:g(h(D_perp,D_theta),JD_perp) Symmetric Warped Product Hemi-Slant Submanifolds}. Similarly,}
			g(h(X,Z),FW)&=g(\overline{\nabla}_XZ,JW)-g(\overline{\nabla}_XZ,PW)\\
			&=-g(J\overline{\nabla}_XZ,W)-g(\nabla_XZ,PW)\\
			&=-g(\overline{\nabla}_XJZ,W)\\
			&=g(\mathfrak{A}_{FZ}X,W)
			\intertext{which implies \eqref{eq:g(h(D_perp,D_theta),FD_theta) Symmetric Warped Product Hemi-Slant Submanifolds}. Repeating the above calculation, we have}
			g(h(X,Z),FW)&=g(\overline{\nabla}_ZX,JW-PW)\\
			&=-g(J\overline{\nabla}_ZX,W)-g(\nabla_ZX,PW)\\
			&=-g(\overline{\nabla}_ZJX,W)-\frac{1}{2}g(JB,X)g(Z,W)-\frac{1}{2}g(B,X)g(JZ,W)\\
			&=g(\mathfrak{A}_{JX}Z,W)-\frac{1}{2}g(JB,X)g(Z,W)-\frac{1}{2}g(B,X)g(PZ,W)
			\intertext{Using \eqref{eq:g(h(D_perp,D_theta),FD_theta) Symmetric Warped Product Hemi-Slant Submanifolds} and comparing symmetric and skew symmetric terms in $Z$ and $W$ we have,}
			g(h(X,Z),FW)&=g(h(Z,W),JX)-\frac{1}{2}g(JB,X)g(Z,W)
			\intertext{which shows \eqref{eq:Identity g(h(D_theta,D_theta),JD_perp) and g(h(D_perp,D_theta),FD_theta) Warped Product Hemi-Slant Submanifolds} and}
			0&=\frac{1}{2}g(B,X)g(PZ,W)
			\intertext{which shows \eqref{eq:Identity g(B,D_perp) Warped Product Hemi-Slant Submanifolds}. Finally,}
			g(h(X,Y),JX_1)&=g(\overline{\nabla}_XY,JX_1)\\
			&=-g(J\overline{\nabla}_XY,X_1)\\
			&=-g(\overline{\nabla}_XJY,X_1)-\frac{1}{2}g(X,X_1)g(JB,Y)+\frac{1}{2}g(X,Y)g(JB,X_1)\\
			&=g(\mathfrak{A}_{JY}X,X_1)-\frac{1}{2}g(X,X_1)g(JB,Y)+\frac{1}{2}g(X,Y)g(JB,X_1)
		\end{align*}
		which gives \eqref{eq:Identity g(h(D_perp,D_perp),FD_perp) Warped Product Hemi-Slant Submanifolds}.
	\end{proof}
	\begin{rem}\label{rem:Local Orthonormal Basis of Warped Product Hemi-Slant Submanifolds}
		Given a warped product hemi-slant submanifold $M_\theta\times\,_\lambda M_\perp$ of an l.c.K manifold $\widetilde{M}^{2n}$, let $\{X_i\}_{i=1}^{p}$ and $\{Z_j,\beta PZ_j\}_{j=1}^{q}$ respectively be local orthonormal frames of $TM_\perp$ and $TM_\theta$. Then a local orthonormal frame of $\widetilde{M}^{2n}$ is
		\[\begin{aligned}
			\left\lbrace \widehat{X_i}=\frac{X_i}{\lambda}\right\rbrace\cup\left\lbrace\widehat{Z_j}=Z_j,\widehat{PZ_j}=\beta PZ_j\right\rbrace\cup\left\lbrace \widehat{JX_i}=\frac{JX_i}{\lambda}\right\rbrace\\\cup\left\lbrace\widehat{FZ_j}=\alpha FZ_j,\widehat{FPZ_j}=\alpha\beta FPZ_j\right\rbrace\cup\left\lbrace\widehat{\xi_k},\widehat{J\xi_k}\right\rbrace
		\end{aligned}\]
		where $\alpha=\csc\theta$, $\beta=\sec\theta$ and
		\begin{align*}
			&\left\{\widehat{X_i}:1\leq i\leq n_1\right\}\text{ is an orthonormal basis of }\mathcal{D}^\perp\\
			&\left\{\widehat{Z_j}, \widehat{PZ_j}:1\leq j\leq n_2\right\}\text{ is an orthonormal basis of }\mathcal{D}^\theta\\
			&\left\{\widehat{JX_i}:1\leq i\leq n_1\right\}\text{ is an orthonormal basis of }J\mathcal{D}^\perp\\
			&\left\{\widehat{FZ_j}, \widehat{FPZ_j}:1\leq j\leq n_2\right\}\text{ is an orthonormal basis of }F\mathcal{D}^\theta\\
			&\left\{\widehat{\xi_k},\widehat{J\xi_k}:1\leq k\leq \frac{n-n_1-2n_2}{2}\right\}\text{ is an orthonormal basis of }\mu
		\end{align*}
		However, while $X_i\in\mathcal{L}(M_\perp)$ we have $\widehat{X_i}\notin\mathcal{L}(M_\perp)$ in general, as $\lambda$ is a function on $M_\theta$.
	\end{rem}
	\noindent We now give a characterisation for warped product hemi-slant submanifolds of l.c.K. manifolds.
	\begin{theo}\label{th:Warped Prod Hemi-Slant Characterisation - M_perp x_lambda M_theta}
		Let $M^m$ be a hemi-slant submanifold of an l.c.K. manifold $\widetilde{M}^{2n}$. Then the following are equivalent 
		\begin{enumerate}
			\item $M^m$ is a warped product hemi-slant submanifold \mbox{$M_\theta\times\,_\lambda M_\perp$} of $\widetilde{M}^{2n}$\label{Warped Prod Hemi-Slant Characterisation}
			\item For every $X,Y\in\mathcal{L}(M_\perp)$ and $Z,W\in\mathcal{L}(M_\theta)$ we have
			\begin{equation}\label{eq:Warped Prod Hemi-Slant Characterisation 1}
				\begin{aligned}
					g(\mathfrak{A}_{JX}Z-\mathfrak{A}_{FZ}X,Y)&=\left(\frac{1}{2}\omega(JZ)-Z(\ln\lambda)\right)g(X,Y)\\
					\omega(\mathcal{D}^\perp)=\{0\}\text{ and }&g(\mathfrak{A}_{FPZ}X,W)+g(\nabla^\perp_WFZ,JX)=0
				\end{aligned}
			\end{equation}
			for some smooth function $\lambda:M_\theta\to(0,\infty)$.\label{Warped Prod Hemi-Slant Characterisation - Condition 1}
			\item For every $X\in\mathcal{L}(M_\perp)$ and $Z\in\mathcal{L}(M_\theta)$ we have
			\begin{equation}\label{eq:Warped Prod Hemi-Slant Characterisation 2}
				\omega(\mathcal{D}^\perp)=\{0\}\text{ and }\nabla_XZ=\nabla_ZX=Z(\ln\lambda)X
			\end{equation}\label{Warped Prod Hemi-Slant Characterisation - Condition 2}
		\end{enumerate}
		Also, in this case we have the mean curvature vector $\mathcal{H}$ of $M_\perp$ in $M^m$ is 
		\begin{equation}\label{eq:Warping Function Warped Prod Hemi-Slant Characterisation}
			\mathcal{H}=-\text{grad}(\ln\lambda)
		\end{equation}
	\end{theo}
	\begin{proof}
		\mbox{\eqref{Warped Prod Hemi-Slant Characterisation}$\Leftrightarrow$\eqref{Warped Prod Hemi-Slant Characterisation - Condition 1}} This follows from Theorem \ref{th:Hemi-Slant D_perp}, Theorem \ref{th:Hemi-Slant D_theta} and the fact that $\text{grad}\ln\lambda\in\mathcal{L}(M_\theta)$ which implies for all $X\in \mathcal{L}(M_\perp)$ and $Z\in \mathcal{L}(M_\theta)$
		\begin{align*}
			g(\nabla_X(\text{grad}\ln\lambda),Z)&=XZ(\ln\lambda)-g(\text{grad}\ln\lambda,\nabla_XZ)\\
			&=[X,Z](\ln\lambda)-\nabla_XZ(\ln\lambda)\hspace{0.4cm}(\text{as }X(\ln\lambda)=0)\\
			&=-\nabla_ZX(\ln\lambda)\\
			&=g(X,\nabla_Z(\text{grad}(\ln\lambda)))\\
			&=0\hspace{0.4cm}
		\end{align*}
		as $\mathcal{D}^\theta$ is totally geodesic.\par\noindent
		\mbox{\eqref{Warped Prod Hemi-Slant Characterisation}$\Leftrightarrow$\eqref{Warped Prod Hemi-Slant Characterisation - Condition 2}} Let \mbox{$M=M_\theta\times\,_\lambda M_\perp$} be a warped product hemi-slant submanifold. Then \eqref{eq:Warped Prod Hemi-Slant Characterisation 2} and \eqref{eq:Warping Function Warped Prod Hemi-Slant Characterisation} follow from \eqref{eq:Connection Warped Product on M_1 x M_2}.\par
		Conversely, let $M^m$ be a hemi-slant submanifold of an l.c.K. manifold $\widetilde{M}^{2n}$ such that \eqref{eq:Warped Prod Hemi-Slant Characterisation 2} holds. Then for all $X,Y\in \mathcal{L}(M_\perp)$ and $Z,W\in \mathcal{L}(M_\theta)$ we have
		\begin{align*}
			g([X,Y],Z)&=g(\nabla_XY-\nabla_YX,Z)\\
			&=-g(\nabla_XZ,Y)+g(\nabla_YZ,X)\\
			&=-Z(\ln\lambda)g(X,Y)+Z(\ln\lambda)g(Y,X)=0
			\intertext{which implies $\mathcal{D}^\perp$ is involutive.}
			g(\nabla_XY,Z)&=-g(\nabla_XZ,Y)\\
			&=-Z(\ln\lambda)g(X,Y)\\
			&=-g(X,Y)g(\text{grad}(\ln\lambda),Z)
			\intertext{which implies leaves of $\mathcal{D}^\perp$ are totally umbilical in $M$ with mean curvature vector $-\text{grad}(\ln\lambda)$.}
			g([Z,W],X)&=g(\nabla_ZW-\nabla_WZ,X)\\
			&=-g(\nabla_ZX,W)+g(\nabla_WX,Z)\\
			&=0
			\intertext{which implies $\mathcal{D}^\theta$ is involutive.}
			g(\nabla_ZW,X)&=-g(\nabla_ZX,W)=0
			\intertext{which implies leaves of $\mathcal{D}^\theta$ are totally geodesic in $M$.}
			g\left(\nabla_X\text{grad}(\ln\lambda),Z\right)&=\text{grad}(\ln\lambda)(\ln\lambda)g(X,Z)=0
		\end{align*}
		which implies $\text{grad}(\ln\lambda)$ is parallel in the normal bundle of $M_\perp$ in $M$.\par
		Hence by Theorem \ref{th:Characterisation - Warped Product}
		we have \mbox{$M=M_\theta\times\,_\lambda M_\perp$} is a warped product hemi-slant submanifold.
	\end{proof}
	We conclude our study of warped product hemi-slant submanifolds of l.c.K. manifolds by giving an inequality for the norm of the second fundamental form.
	\begin{theo}\label{th:Length h - Warped Product Hemi-Slant}
		Let \mbox{$M= M_\perp\times \,_\lambda M_\theta$} be a warped product hemi-slant submanifold in an lcK manifold $(\widetilde{M}^{2n},J,g )$. Then the norm of the second fundamental form satisfies the inequality
		\begin{align}\label{eq:Length h - Warped Product Hemi-Slant}
			||h||^2\geq\:&
			\frac{n_1+n_2-1}{2}\|B|_{J\mathcal{D}^\perp}\|^2+2g\left(H_{\mathcal{D}^\perp}|_{J\mathcal{D}^\perp},B|_{J\mathcal{D}^\perp}\right)+g\left(H_{\mathcal{D}^\theta}|_{J\mathcal{D}^\perp},B|_{J\mathcal{D}^\perp}\right)-K
		\end{align}
		where $n_1=\dim_{\mathbb{R}}\mathcal{D}^\perp$, $2n_2=\dim_{\mathbb{R}}\mathcal{D}^\theta$, $H_{\mathcal{D}^\perp}$ and $H_{\mathcal{D}^\theta}$ are respectively the components of the mean curvature vector $H$ of $M$ in $\widetilde{M}^{2n}$ along $\mathcal{D}^\perp$ and $\mathcal{D}^\theta$ and given any orthonormal basis $\left\{\widehat{X_i}:1\leq i\leq n_1\right\}$ of $\mathcal{D}^\perp$, \mbox{$K=2\sum_{i}g\left( h\left(\widehat{X_i},\widehat{X_i}\right),\widehat{JX_i}\right)g\left(B,\widehat{JX_i}\right)$}  .\par
		\noindent If equality holds then we have
		\begin{itemize}
			\item $\text{Image}(h)\subseteq(J\mathcal{D}^\perp\oplus F\mathcal{D}^\theta)$, and
			\item $M_\theta$ is totally geodesic in $\widetilde{M}^{2n}$, if and only if, $M$ is mixed-totally geodesic in $\widetilde{M}^{2n}$.
		\end{itemize}
	\end{theo}
	The proof follows on the same lines as that of Theorem \ref{th:Length h - Hemi-Slant Warped Product}.
	\par\bigskip \noindent
	\textbf{Acknowledgement:} The authors would like to thank the referee(s) for their invaluable criticism and suggestions towards improving the paper.

\end{document}